\documentclass{amsart}
\usepackage{amsmath}
\usepackage{enumitem}
\usepackage{amscd}
\usepackage{amsfonts}
\usepackage[small,nohug]{diagrams}
\usepackage{bm}
\usepackage{pstricks}
\usepackage{pst-plot}

\newtheorem{theorem}{Theorem}[section]

\newtheorem{proposition}[theorem]{Proposition}

\theoremstyle{remark}

\theoremstyle{definition}

\numberwithin{equation}{section}

\DeclareMathOperator{\R}{R}

\DeclareMathOperator{\GP}{GP}
\DeclareMathOperator{\dist}{d}

\begin{document}

\title[Uniformly joinable, locally uniformly joinable, and weakly chained]{Uniformly joinable, locally uniformly joinable, and weakly chained uniform spaces}
\author{Brendon LaBuz}
\address{Saint Francis University, Loretto, PA 15940}
\email{blabuz@@francis.edu}
\date{\today}

\begin{abstract}
Brodskiy, Dydak, LaBuz, and Mitra introduced the concepts of uniform joinability and local uniform joinability for uniform spaces when developing their theory of generalized uniform covering maps which was motivated by a paper by Berestovskii and Plaut. (Local) uniform joinability can be thought of as analogous to (local) path connectedness. A chain connected locally uniformly joinable uniform space is uniformly joinable. This note gives an example of a metric space that is uniformly joinable but not locally uniformly joinable. 

Plaut recently defined the concept of a weakly chained uniform space. We show that a weakly chained metrizable uniform space is locally uniformly joinable. Since local uniform joinability is equivalent to pointed 1-movability for metric continua, we find that weak chainability is equivalent to pointed 1-movability for such spaces.
\end{abstract}

\maketitle
\tableofcontents

\section{Introduction}
The concepts of uniform joinability and local uniform joinability were defined in \cite{Rips} in the context of creating generalized uniform covering maps. Traditional covering space theory does not apply to locally complicated spaces. For such spaces one can consider chains instead of paths, and the inverse limit of finer and finer chains gives generalized paths. These chains require a way of measuring closeness across a space so the natural venue for such a construction is the category of uniform spaces.

Given a uniform space $X$ and an entourage $E$ of $X$, an $E$-chain in $X$ is a sequence of points $x_1,\ldots, x_n$ in $X$ such that $(x_i,x_{i+1})\in E$ for each $i<n$. We define $E$-homotopies of $E$-chains by appealing to Rips complexes. Given an entourage $E$ of $X$ the Rips complex $\R(X,E)$ is the subcomplex of the full complex over $X$ whose simplices are finite $E$-bounded subsets of $X$. A subset $A$ of $X$ is called $E$-bounded if for each $x,y\in A$, $(x,y)\in E$. Equivalently, $A\times A\subset E$. Any $E$-chain $x_1,\ldots ,x_n$ determines a path in $R(X,E)$ by simply joining successive terms $x_i,x_{i+1}$ by an edge path, a path along the edge joining $x_i$ and $x_{i+1}$. Then $E$-homotopies between $E$-chains can be defined in terms of homotopies between paths in $\R(X,E)$.\footnote{Since the identity function $K_w\to K_m$, $K$ a simplicial complex, from $K$ equipped with the CW (weak) topology to $K$ equipped with the metric topology is a homotopy equivalence, it does not matter which topology is chosen for $\R(X,E)$.} Two $E$-chains starting at the same point and ending at the same point are said to be $E$-homotopic if the corresponding paths in $\R(X,E)$ are homotopic relative endpoints. This homotopy in $\R(X,E)$ can be chosen to be simplicial. Thus, two
$E$-chains are $E$-homotopic if one can move from one to the other by a finite sequence of the addition or deletion of a point in the $E$-chain (while keeping the endpoints fixed). Given an $E$-chain $c$, denote the $E$-homotopy class of $c$ as $[c]_E$.

A generalized path is a collection of homotopy classes of chains $\alpha=\{[c_E]_E\}_E$ where $E$ runs over all entourages of $X$ and for any $F\subset E$, $c_F$ is $E$-homotopic to $c_E$. We write $\alpha_E=[c_E]_E$. Inverses and concatenations of generalized paths are defined in the obvious way. The set of generalized paths in $X$ is denoted as $\GP(X)$. We also have the pointed version $\GP(X,x_0)$ of all generalized paths in $X$ starting at some basepoint $x_0$. Given an entourage $E$ of $X$, define $E^*$ to be the set of all pairs $(\alpha,\beta)$ of elements of $GP(X,x_0)$ such that $\alpha_E^{-1}*\beta_E$ is $E$-homotopic to the $E$-chain $x,y$ where $x$ is the endpoint of $\alpha$ and $y$ is the endpoint of $\beta$. Such a generalized path is called $E$-short. The set of all $E^*$ as $E$ runs over all entourages of $X$ forms a basis for a uniform structure on $\GP(X,x_0)$. 

A uniform space $X$ is called uniformly joinable if every pair of points in $X$ can be joined by a generalized path. Equivalently, the endpoint map $\GP(X,x_0)\to X$ is surjective. Note any path in $X$ induces a generalized path so if $X$ is path connected then it is uniformly joinable. A uniform space $X$ is called locally uniformly joinable if for each entourage $E$ of $X$ there is an entourage $F$ of $X$ such that if $(x,y)\in F$ then $x$ and $y$ can be joined by an $E$-short generalized path.

A metric space $X$ induces a uniform structure as follows. Given $\epsilon>0$, $E_\epsilon$ is defined to be the set of pairs $(x,y)$ of points in $X$ such that $\dist(x,y)\leq \epsilon$. Then the set of all $E_\epsilon$ form a basis for the uniform structure on $X$.

For metric continua, by a theorem of Krasinkiewicz and Minc, uniform joinability and local uniform joinability are equivalent (see \cite[Corollary 6.6]{Rips}). Any chain connected locally uniformly joinable uniform space is uniformly joinaable \cite[Proposition 4.3]{Rips} and it is easy to see that there are spaces that are locally uniformly joinable but not uniformly joinable (consider parallel lines in the plane). We give an example of a metric space that is uniformly joinable but not locally uniformly joinable. It is a path connected subset of the plane with the uniform structure induced by the usual metric.

\section{The Texas circle}

\begin{figure}\label{asymptote}
\psset{unit=7.5mm}
\begin{pspicture*}(3.14156,0)(19.5,1.31831)
\psplot[algebraic,plotpoints=500]{3.1416}{18.5}{sin(x)^2+1/x}
\psline(3.14156,0)(18.5,0)
\psset{dotscale=0.5}
\psdots(18.75,0.5)(19,0.5)(19.25,0.5)
\psset{linewidth=2pt}
\psline(3.14156,0)(3.14156,0.31831)
\end{pspicture*}
\caption{The Texas circle}
\end{figure}
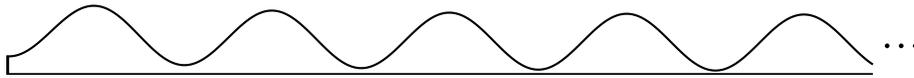

Let $T$ be the union of the graph of $f(x)=\sin^2x+1/x$ over $[\pi,\infty)$, the $x$-axis from $\pi$ to $\infty$, and the vertical line segment from $(\pi,0)$ to $(\pi,1/\pi)$ (see Figure \ref{asymptote}). We call $T$ the Texas circle as it can be considered a large scale version of the Warsaw circle. Give $T$ the uniform structure induced by the usual metric on the plane. Then $T$ is path connected so it is uniformly joinable. Let us see that $T$ is not locally uniformly joinable. Suppose to the contrary that $T$ is locally uniformly joinable. Then there is an entourage $F$ so that if $(x,y)\in F$, $x$ and $y$ can be joined by an $E_{1/2}$-short generalized path. Let $n\in\mathbb N$ be such that if $\dist(x,y)\leq 1/n\pi$ then $(x,y)\in F$. Consider the points $x=(n\pi,1/n\pi)$ and $y=(n\pi,0)$. Then $(x,y)\in F$ so there is an $E_{1/2}$-short generalized path joining them. But the only generalized path joining them is the one induced by the path that travels along the graph of $f$ all the way to the vertical line segment and then back along the $x$-axis and this generalized path is not $E_{1/2}$-short (notice a square of side length $1/2$ fits under the first crest).

To see that there are no other generalized paths joining $x$ and $y$, consider that for $m>n$ and $\epsilon\leq 1/m\pi$, the only $E_\epsilon$-chains joining $x$ and $y$ are ones that travel along the graph of $f$ all the way to the vertical line segment and then back along the $x$-axis, and the ones that travel to the right of $x$ to at least the point $(m\pi,1/m\pi)$ and then hop down to the $x$-axis and back to $y$ (or some combination of these two types of $E_\epsilon$-chains). Now for these types of $E_\epsilon$-chains to form a generalized path we need that for $m'>m$, the $E_{1/m'\pi}$-chain is $E_{1/m\pi}$-homotopic to the $E_{1/m\pi}$-chain. But that cannot always be the case. For fix an $E_{1/m\pi}$-chain $c$. This chain has a right-most point, say with $x$-value $x_m$. Choose $m'$ so that $(m'-1)\pi>x_m$. Then any $E_{1/m'\pi}$-chain joining $x$ to $y$ either travels through the vertical line segment or must travel to at least $(m'\pi,1/m'\pi)$ before hopping down to the $x$-axis and these types of $E_{1/m'\pi}$-chains are not $E_{1/m\pi}$-homotopic to $c$.

\section{Weakly chained spaces}

A uniform space is called chain connected if for every entourage $E$, any pair of points can be joined by an $E$-chain. In the preprint of Plaut \cite{P2}, a uniform space is defined to be weakly chained if it is chain connected and for any entourage $E$ of $X$ there exists an entourage $F\subset E$ such that for any $(x,y)\in F$, there exists arbitrarily fine chains $c$ joining $x$ and $y$ such that $[c]_E=[x,y]_E$ (we can say that the chain is $E$-short).
Clearly a chain connected locally uniformly joinable space is weakly chained. We will see that the converse holds for metrizable spaces.

\begin{proposition}
Let $X$ be a metrizable uniform space. If $X$ is weakly chained then it is locally uniformly joinable.
\end{proposition}

\begin{proof}
Let a metric be given that induces the uniform structure on $X$. For $n>0$, recall $E_{1/n}$ is the entourage consisting of pairs of points of $X$ that are distance at most $1/n$ apart. Let $E$ be an entourage of $X$. We will find an entourage $F$ of $X$ so that if $(x,y)\in F$, $x$ and $y$ can be joined by an $E$-short generalized path. Let $F_1\subset E$ be an entourage such that if $(x,y)\in F_1$, $x$ and $y$ can be joined by arbitrarily fine chains that are $E$-short. Let $F_{2}\subset E_{1/2}$ be an entourage such that if $(x,y)\in F_2$, $x$ and $y$ can be joined by arbitrarily fine chains that are $F_{1}$-short. Inductively choose, for each $n\geq 3$, an entourage $F_{n}\subset E_{1/n}$ such that if $(x,y)\in F_n$, $x$ and $y$ can be joined by arbitrarily fine chains that are $F_{n-1}$-short.

Suppose $(x,y)\in F_2$. We will inductively define $F_n$-chains from $x$ to $y$ that will create a generalized path that is $E$-short. Let $c_1$ be an $F_3$-chain joining $x$ and $y$ that is $F_1$-short. We define $[c_1]_{F_1}$ to be the $F_1$-level of the generalized path we are constructing. Note $[c_1]_{F_1}=[x,y]_{F_1}$ so $c_1$ is $E$-short. 

Now join pairs of successive points of $c_1$ by $F_4$-chains that are $F_2$-short. Let $c_2$ be the $F_4$-chain resulting from concatenating those chains. We define $[c_2]_{F_2}$ to be the $F_2$-level of the generalized path we are constructing. Note $[c_2]_{F_2}=[c_1]_{F_2}$. Using induction, for each $n\geq 3$, join pairs of successive points of $c_{n-1}$ by $F_{n+2}$-chains that are $F_n$-short. Let $c_n$ be the $F_{n+2}$-chain resulting from concatenating those chains. We define $[c_n]_{F_n}$ to be the $F_n$-level of the generalized path we are constructing. Note $[c_n]_{F_n}=[c_{n-1}]_{F_n}$.

These equivalence classes induce a generalized path from $x$ to $y$ that is $E$-short.
\end{proof}

Recall that a chain connected locally uniformly joinable uniform space is uniformly joinable (\cite[Proposition 4.3]{Rips}). Thus a weakly chained uniform space is uniformly joinable. Also, according to \cite[Corollary 6.6]{Rips}, for metric continua, the conditions of uniformly joinable, locally uniformly joinable, and pointed 1-movable are equivalent. Since metric continua are chain connected we can add weakly chained to the list  of equivalent conditions.

\end{document}